\documentclass[reqno,11pt,centertags]{amsart}
\usepackage{amsmath,amsthm,amscd,amssymb,latexsym,upref}
\date{\today}

\newcommand{\bbN}{{\mathbb{N}}}
\newcommand{\bbR}{{\mathbb{R}}}

\newcommand{\bbC}{{\mathbb{C}}}

\renewcommand{\Re}{\text{\rm Re}}

\def\sgn{{\mathrm{sgn}\,}}
\renewcommand{\Im}{\text{\rm Im}}

\allowdisplaybreaks
\numberwithin{equation}{section}
\newtheorem{theorem}{Theorem}[section]

\newtheorem{lemma}[theorem]{Lemma}
\newtheorem{corollary}[theorem]{Corollary}

\theoremstyle{definition}

\newtheorem{remark}[theorem]{Remark}
\newtheorem{conjecture}[theorem]{Conjecture}

\newtheorem{problem}[theorem]{Problem}

\begin{document}

\title[]{Asymptotics of the best  polynomial  approximation of $|x|^p$  and of the best Laurent polynomial
approximation of $\sgn(x)$ on  two symmetric intervals}
\author[]{F. Nazarov, F. Peherstorfer, A. Volberg  and P. Yuditskii}

\thanks{
Partially supported by NSF grant DMS-0501067  the Austrian Founds
FWF, project number: P16390--N04 and {\it Marie Curie International
Fellowship} within the 6-th European Community Framework Programme,
 Contract MIF1-CT-2005-006966}

\date{\today}


\maketitle

\begin{abstract} We present a new method that allows us to get a direct proof
of the classical Bernstein asymptotics for the error of the best
uniform polynomial approximation of $|x|^p$ on two symmetric
intervals. Note, that in addition, we get asymptotics  for the
polynomials themselves under a certain renormalization. Also, we
solve  a problem on asymptotics of the best approximation of
$\sgn(x)$ on  $[-1,-a]\cup[a,1]$ by Laurent polynomials.

\end{abstract}

\section{Introduction}%

By S. N. Bernstein  \cite{Bern3}, \cite[Ch. II, Sect. 5]{Bernstb},
see also \cite{Akha}, it follows
 that the error $E_n(p,a)$ of the best uniform approximation
of $|x|^p$, $p$  not an even integer, on $[-1,-a]\cup[a,1]$ by
polynomials of degree $n=2m$ the following limit exists:
\begin{equation}\label{f1}
\lim_{m\to \infty}\left(\frac{1+a}{1-a}\right)^{m+1}
m^{\frac{p}{2}+1}E_n(p,a)={a}^{\frac{p}{2}-1}
\frac{(1+a)^2}{2\left|\Gamma\left(-\frac p 2\right)\right|}.
\end{equation}
Indeed, see Appendix 2, \eqref{f1} may be derived easily from the
error of the best approximation of $1/(b+x)^s$, $s\not=0$. We
mention that due to Chebyshev, see \cite[chap. II, No. 37]{Akha} or
\cite[p. 120]{Bern3}, for $s=1$ even the best polynomial
approximation is known explicitly.

Note that for $a=0$, that is for  the approximation of $|x|^p$ on
$[-1,1]$, the asymptotics for $ E_n(p)= E_n(p,0)$ is still the open
famous Bernstein
 problem, though the existence  of the limit
 \begin{equation}\label{E}
    \lim_{n\to\infty}n^p E_n(p)=\mu(p)>0
 \end{equation}
was shown in \cite{Bern1,Bern2}. In particular, one can not just put
$a=0$ in \eqref{f1}, even the growth (exponent for $n=2m$) is
different. For recent progress concerning the Bernstein problem, see
D. Lubinsky \cite{L, L2}. A survey on Bernstein constant theorems is
given in \cite{G1}, see also   \cite{G2}.

 For results of the type \eqref{E}, where $|x|^p$ (resp.
 $|x-x_0|^p$) is approximated on a system of several intervals, see
 \cite{V, T}, \cite[Sect. 10]{T2}.

The  interest to this remarkable problem (see, e.g. \cite{L}) was
boosted by the  recent result of H. Stahl \cite{Sta}, who completed
a long line of studying of the analogous problem for uniform {\it
rational} approximation of $|x|^p$ on $[-1,1]$  with the remarkable
explicit answer:
$$\lim_{n\to\infty}\exp(\pi\sqrt{p n})E^r_n(p)=
2^{2+p}|\sin(\pi p/2)|,$$ where $E^r_n(p)$ is the error of the best
rational approximation.

E. I. Zolotarev \cite{Zol,Akh} found an explicit expression, in
terms of elliptic functions, of the rational function of given
degree which is uniformly closest to $\sgn(x)$ on
$[-1,-a]\cup[a,1]$.

 We call a rational function of the form
$$
f(x)=\frac{a_{-l}}{x^{l}}+...+{a_{n}}{x^{n}}
$$
a Laurent polynomial of degree $(l,n)$.

\begin{problem}\label{p1.1} For $k,m\in\bbN$,
find the best
approximation of the function $\sgn(x)$,
 $|x|\in [a,1]$, by Laurent polynomials of degree $(2k-1,2m-1)$
and   the
approximation error $L^k_m(a)$.
\end{problem}

\begin{remark}\label{r1}
Problem \ref{p1.1},   in a trivial way, is related to the following
 weighted polynomial approximation problem:
\begin{equation}\label{221}
E^*_n(p,a)=L_m^k(a)=\inf_{\{P:\deg P\le
2(m+k-1)\}}\sup_{|x|\in[a,1]}
\left\vert\frac{|x|^{2k-1}-P(x)}{x^{2k-1}}\right\vert,
\end{equation}
where $a\in (0,1)$, $k,m\in \bbN$ and $p=2k-1$, $n=2(k+m-1)$. Also,
it is trivial that the extremal function $f=f(x;k,m;a)$ is odd
 in Problem \ref{p1.1} and the
extremal polynomial in  \eqref{221} is even.
\end{remark}

\begin{problem}\label{p1.3} For an integer $m$ and a real $p$,
$2m>p> 0$, find the polynomial of degree $2m$ of the best
approximation to $|x|^p$,
 $|x|$ belongs to $[a,1]$.
\end{problem}

In \cite{EYU} A. Eremenko and the last author solved the standard
polynomial case and the Laurent case of this problem was considered
in \cite{PYU}. In particular, Lemma \ref{l2.1}, Theorems \ref{t3.1}
and \ref{t5.1}, and a weak version of Theorem \ref{t4.1} were proved
in preprint \cite{PYU}. However, in both works, in the last step we
were enforced to use Bernstein's result similar to \eqref{f1}. Here
we, finally, close this approach, and, as it was mentioned before,
obtain in a natural way asymptotics for the approximation error and
simultaneously for the extremal function under a certain
renormalization.

The main steps of the method are (with respect to  Problem
\ref{p1.1}):
\begin{itemize}
\item[1.] For each particular $k$ and $m$ we reveal the structure of the
extremal function by representing it with the help of an explicitly
given conformal mapping.
\item[2.] The system of conformal mappings ($k$ is fixed, $m$ is a
parameter) converges (in the Caratheodory sense) after an
appropriate renormalization. The limit map does not depend on $a$,
thus we obtain asymptotics for $L^k_m(a)$ in terms of $a$--depending
parameters, that we use for renormalization, (an explicit formula)
and a $k$--depending constant, say $Y_k$, which is a certain
characteristic of the final conformal map (kind of capacity).
\item[3.] Using a special representation for bounded Nevanlinna
functions we get the explicit formula for the final conformal
mapping, in particular, for the constant $Y_k$.
\end{itemize}

With a slight modification we apply the method to asymptotics
related to Problem \ref{p1.3}, see Sections 5 and 6. Note, that up
to the last step in this problem, we can follow  the same program in
the most intriguing case $a=0$. However, on the contrary to the
linear equation appearing in the considered case $a>0$, see
\eqref{lineq2} (or \eqref{lineq} for Problem \ref{p1.1}), we get a
kind of quadratic equation \eqref{al4} involving an unknown
function, its Hilbert transform and an independent variable. The
trigonometric form \eqref{trigeq} might be preferable for the
equation. But, in any case, at the moment we are unable to find a
way to get its explicit solution and for this reason we do not
discuss this subject in the main part of the paper and just
formulate corresponding conjecture in Appendix 1.

\bigskip

\noindent{\bf Acknowledgment.} We are thankful to Alex Eremenko for
friendly conversations during the writing of this paper and referees
for remarks that help to improve presentation and to add essential
references related to the subject.

\section{Special Conformal Mappings}

In this section we introduce certain special conformal mappings that we
need in what follows. They are marked by a natural parameter $k$, but in
this section $k$ can be just real, $k>1/2$.

For given  $k$, consider the domain
\begin{equation}\label{f0}
\Pi_k=\bbC_+\setminus\{w: \Re w=-\log t,\ |\Im w-k\pi|\le \arccos t,
\ t\in (0,1]\}
\end{equation}
Define the conformal map
$$
H_k:\bbC_+\to \Pi_k
$$
normalized by $H_k(0)=\infty_1$,
$H_k(\infty)=\infty_2$ (on the boundary we have two infinite points
that we denote respectively $\infty_1,\infty_2$), and moreover
$$
H_k(\zeta)=\zeta+..., \quad \zeta\to \infty,
$$
(that is the leading coefficient is fixed).
By $D_k$ we denote the positive number such that $H_k(-D_k)=0$.

Recall that a Nevanlinna function $G(z)$ ($\Im\, G(z)>0$, for $\Im\,
z>0$)  possesses the integral representation, see e.g. \cite{KN},
\begin{equation}\label{nf}
G(z)=A z+B+\int\left(\frac{1}{t-z}-\frac{t}{1+t^2}\right)d\sigma(t),
\end{equation}
where $A>0$, $B\in \bbR$, $\sigma$ is a positive measure on the real
axis such that $\int\frac{d\sigma(t)}{1+t^2}<\infty$. Moreover
\begin{equation}\label{nfpar}
    A=\lim_{z=iy, y\to\infty}\frac{G(z)}z,\quad
\sigma(x_2)-\sigma(x_1)=\lim_{\epsilon \to 0} \frac 1
\pi\int_{x_1}^{x_2}\Im\, G(x+i\epsilon)\,dx.
\end{equation}

Therefore  for $H_k$ we have the following integral representation
\begin{equation}
H_k(\zeta)=\zeta+D_k+\int_0^\infty
\left(\frac{1}{t-\zeta}-\frac{1}{t+D_k}\right)\rho_k(t) dt,
\end{equation}
where $\rho_k(t)=\frac{1}{\pi}\Im H_k(t)$. Evidently
$\rho_k(t)\to k+\frac 1 2$, $t\to +\infty$.
\begin{lemma} \label{l2.1} The function $H_k$ possesses the asymptotic
\begin{equation}\label{2.3}
\lim_{\zeta\to-\infty}
\left\{
H_k(\zeta)-\zeta+\left(k+\frac 1 2\right)\log(-\zeta)
\right\}=Y_k,
\end{equation}
where
\begin{equation}\label{new2.3}
Y_k:=D_k+ \left(k+\frac 1 2\right)\log D_k -\int_0^\infty
\frac{\rho_k(t)-\left(k+\frac 1 2\right)}{t+D_k} dt.
\end{equation}

\end{lemma}
\begin{proof}
Since
\begin{equation}
\int_0^\infty
\left(\frac{1}{t-\zeta}-\frac{1}{t+D_k}\right)
\left(\rho_k(t)-\left(k+\frac 1 2\right)\right)
dt\to
-\int_0^\infty
\frac{\rho_k(t)-\left(k+\frac 1 2\right)}{t+D_k}
dt,
\end{equation}
as $\zeta\to-\infty$ and
\begin{equation}
\left(k+\frac 1 2\right)\int_0^\infty
\left(\frac{1}{t-\zeta}-\frac{1}{t+D_k}\right)
dt=-\left(k+\frac 1 2\right)(\log(-\zeta)-\log D_k)
\end{equation}
we get \eqref{2.3}.
\end{proof}

As it was mentioned in the Introduction (step 3 of our method), we
will use a certain special representation for a Nevanlinna function
$H_k(\zeta)$. We note that a Nevanlinna function $F(z)$ with the
imaginary part in $[0,\pi]$ has the form $ F(z)=\log G(z)$, where
$G(z)$ is another nontrivial Nevanlinna function, see e.g.
\cite{KN}.
  Based on this remark we
get the following corollary of the previous lemma.
\begin{corollary}
The function $H_k$ possesses the representation
\begin{equation}\label{f2}
 H_k(\zeta)=\zeta-\left(k-\frac 1 2\right)\log(-\zeta)+ \log\left\{
 \frac 1 \pi
\int_0^\infty\frac{\tau_k(t)\,dt}{t-\zeta}\right\}.
\end{equation}
\end{corollary}
\begin{proof} As it follows from \eqref{new2.3} and \eqref{f0}
$$
\Im\left\{H_k(\zeta)-\left(\zeta-\left(k-\frac 1
2\right)\log(-\zeta)\right)\right\}\in [0,\pi],\quad \Im\,\zeta>0.$$
Using the representation \eqref{nf} and \eqref{nfpar} we get
\eqref{f2} with
$$
\tau_k(x)=\arg\left\{H_k(x)-\left(x-\left(k-\frac 1
2\right)\log(-x)\right)\right\}.
$$
\end{proof}

\begin{theorem} The function $H_k(z)$ is of the form
\begin{equation}\label{f3}
 H_k(\zeta)=\zeta-\left(k-\frac 1 2\right)\log(-\zeta)+ \log\left\{\frac 1 \pi
\int_0^\infty\frac{t^{k-\frac 1 2}e^{-t}\,dt}{t-\zeta}\right\}.
\end{equation}
in particular,
\begin{equation}\label{f4}
    Y_k=\log{\Gamma\left(k+\frac 1 2 \right)} -\log\pi.
\end{equation}

\end{theorem}

\begin{proof}

We note that for $\zeta=\xi+i\eta$ the curve in \eqref{f0} is given
by the equation
$$
\Re \{e^{w-i\pi k}\}|_{\zeta=\xi+i0}=1,\quad \xi>0.
$$
We use here the representation \eqref{f2} and, thus, get
\begin{equation}\label{lineq}
    \Im \left .\left\{e^\xi \xi^{-(k-\frac 1 2)} \frac 1 \pi
\int_0^\infty\frac{\tau_k(t)\,dt}{t-\zeta}\right\}\right|_{\zeta=\xi+i0}=1,\quad
\xi>0.
\end{equation}
Therefore \eqref{f3} is proved.

Now we use the standard asymptotic formula
$$
\frac 1{\pi}\int_0^\infty\frac{\tau_k(t)\,dt}{t-\zeta}=\frac{\frac 1
\pi\int_0^\infty \tau_k(t)\,dt}{-\zeta}+ \frac{\frac 1
\pi\int_0^\infty t\tau_k(t)\,dt}{-\zeta^2}+...
$$
Therefore
\begin{equation}\label{newf3}
\begin{split}
H_k(\zeta)=&\zeta-\left(k-\frac 1 2\right)\log(-\zeta)+
\log\left\{\frac{\frac 1 \pi
\int_0^\infty t^{k-\frac 1 2}e^{-t}\,dt}{-\zeta}(1+o(1))\right\}\\
=&\zeta-\left(k+\frac 1 2\right)\log(-\zeta)+ \log\left\{\frac 1 \pi
\int_0^\infty t^{k-\frac 1 2}e^{-t}\,dt\right\}+o(1).
\end{split}
\end{equation}
Due to \eqref{2.3} we get \eqref{f4}.

\end{proof}

\begin{remark} The conformal map on the domain
\begin{equation}\label{f5}
\Pi_0=\bbC_+\setminus\{w: \Re w=-\log t,\ \Im w\le \arccos t, \ t\in
(0,1]\}
\end{equation}
is important for a description of the standard polynomial
approximation of $\sgn(x)$ \cite{EYU}. Though it requires a special
consideration, the formal extension of \eqref{f3} to the case $k=0$,
that is, the formula
\begin{equation}\label{f6}
 H_0(\zeta)=\zeta+\frac 1 2\log(-\zeta)+ \log\left\{\frac 1 \pi
\int_0^\infty\frac{t^{-\frac 1 2}e^{-t}\,dt}{t-\zeta}\right\}
\end{equation}
holds true for the map normalized by $H_0(0)=0$.
\end{remark}

\section{Extremal Problem}

For a parameter $B>0$ and $k,m\in \bbN$, $\Omega_m^k(B)$ denotes the
subdomain of the half strip
$$
\{w=u+iv: v>0, \ 0<u<(k+m)\pi\}
$$
that we obtain by deleting the subregion
\begin{equation}\label{8aug}
\{w=u+iv: |u-\pi k|\le\arccos\left(\frac{\cosh B}{\cosh v}\right),
\ v\ge B\}.
\end{equation}
Let $\phi(z)=\phi(z;k,m;B)$ be the conformal map of the first quadrant
onto $\Omega_m^k(B)$ such that
$\phi(0)=\infty_1$, $\phi(1)=(k+m)\pi$, $\phi(\infty)=\infty_2$.
Let $a=\phi^{-1}(0)$. Then $a$ is a continuous strictly increasing
function of $B$, moreover $\lim_{B\to 0} a(B)=0$ and
$\lim_{B\to\infty}a(B)=1$. Thus we may consider the inverse function
$B(a)=B^k_m(a)$, $a\in(0,1)$.
\begin{theorem} \label{t3.1} The error of the best approximation in Problem \ref{p1.1} is
\begin{equation}\label{3.2}
L_m^k(a)=\frac{1}{\cosh B^k_m(a)}
\end{equation}
and the extremal function is of the form
\begin{equation}\label{f10}
    f(x;k,m;a)=1-(-1)^k L_m^k(a)\cos\phi(x;k,m;B(a)), \quad x>0.
\end{equation}
\end{theorem}

\begin{proof}
 By inspection of the boundary correspondence, we conclude that
$f=1-(-1)^k L\cos\phi$ is real on the positive ray and pure
imaginary on the positive imaginary ray. So by two reflections $f$
extends to a function analytic in $\Bbb C\setminus\{0\}$. The
extended function evidently satisfies
$$\overline{f(\overline{z})}=
f(z)\quad\mbox{and}\quad-\overline{f(-\overline{z})}=f(z),$$ so we
conclude that $f$ is odd. The region  $\Omega_m^k(B)$ is close to
the strip
$$
\left\{ w:\Re w\in \left(0,\pi\left(k-\frac 1
2\right)\right)\right\}
$$
as $\Im
w\to\infty_1$, and to the strip
$$
\left\{ w:\Re w\in \left(\pi\left(k+\frac 1
2\right),\pi(k+m)\right)\right\}
$$
as $\Im w\to\infty_2$. So $\phi\sim (2k-1)\log 1/z,\; z\to 0$,
$\phi\sim (2m-1)\log z,\; z\to \infty$, and, therefore, $f$ is a
Laurent polynomial of degree $(2k-1,2m-1)$.  Now we note that the
graph of $f$ alternates $k+m+1$ times on $[a,1]$ between $1-L$ and
$1+L$. That a Laurent polynomial with such graph is the unique
extremal for Problem \ref{p1.1} follows from the general theorem of
Chebyshev on the uniform approximation of continuous functions
\cite[Ch. II]{Akha}.

Finally, we have to note that on the imaginary axis the extremal
function has precisely one zero (there are no critical points and
the behavior at $i0$ and at $i\infty$ is evident). At this point
$\phi=k\pi+iB$ and we have \eqref{3.2}.
\end{proof}

\section{Asymptotics}
\begin{theorem} \label{t4.1} The following limit exists
\begin{equation}\label{4.1}
\begin{split}
\lim_{m\to\infty}\left\{B^k_m(a)-\left(m-\frac 1 2\right)
\log\frac{1+a}{1-a}-\left(k+\frac 1 2\right)\log (2m-1)
\right\}\\=\left(k+\frac 1
2\right)\log{\frac{a}{1-a^2}}-\log\frac{\Gamma(k+1/2)}{\pi}.
\end{split}
\end{equation}
Moreover, uniformly on compact subsets of the positive half--axis,
\begin{equation}\label{f12}
\begin{split}
\lim_{m\to\infty}&
    f\left(\sqrt{\frac{2a}{2m-1}}\lambda;k,m;a\right)\\
    =& 1+\frac {(-1)^{k+1}} \pi \int_0^\infty
    \left(\frac{\mu}{\lambda}\right)^{2k-1}e^{-(\lambda^2+\mu^2)}\frac
    {2\mu d\mu}
    {\lambda^2+\mu^2}.
    \end{split}
\end{equation}
\end{theorem}

\begin{proof}
We use the symmetry principle and make a convenient changes of
variable to have a conformal map $\Phi_m(Z)=\Phi(Z;k,m;B)$ of the
upper $Z$--plane
\begin{equation}\label{f.7}
    Z=C_m\sqrt{\frac{z^2-a^2}{z^2-1}}
\end{equation}
in  the region
$$
i(\Omega_m^k(B)\cup\overline{\Omega_m^k(B)})\cup(0,i\pi(m+k)).
$$
This conformal map has the following boundary correspondence
$$
\Phi_m:(-C_m,-A_m,0,A_m,C_m)\to
(-\infty_2,-\infty_1,0,\infty_1,\infty_2),
$$
here $A_m=aC_m$ and the parameter $C_m$ will be chosen a bit later.

For $\Phi_m$ we have the following integral representation
$$
\Phi_m(Z)=\left(m-\frac 1 2\right)\log\frac{1+\frac{Z}{C_m}}
{1-\frac{Z}{C_m}}+\int_{A_m}^\infty
\left[\frac{1}{X-Z}-\frac{1}{X+Z}\right]v_m(X)\,dX,
$$
where
\begin{equation}
v_m(X)=
\begin{cases}
\frac 1{\pi}\Im \Phi_m(X),& A_m\le X\le C_m \\
k+\frac 1 2,& X> C_m
\end{cases}
\end{equation}
Put now
$$
H_m^k(\zeta)=\Phi_m(Z)-B_m,\quad Z=A_m+\zeta,
$$
then
\begin{equation*}
    \begin{split}
    H^k_m(\zeta)=&\left(m-\frac 1
2\right)\log\frac{1+a+\frac{\zeta}{C_m}}
{1-a-\frac{\zeta}{C_m}}+\int_{0}^\infty
\left[\frac{1}{t-\zeta}-\frac{1}{t+2 A_m+\zeta}\right]\hat
v_m(t)\,dt \\
-&B_m,
\end{split}
\end{equation*}
where $\hat v_m(t)=v_m(t+A_m)$. Let us rewrite $H^k_m$ in the form that
is close to the integral representation of $H_k$:
\begin{equation}\label{4.2}
\begin{split}
H^k_m(\zeta)=&\left(m-\frac 1 2\right)\log\frac{1+\frac{\zeta}{C_m(1+a)}}
{1-\frac{\zeta}{C_m(1-a)}}+D_k+
\int_{0}^\infty
\left[\frac{1}{t-\zeta}-\frac{1}{t+D_k}\right]\hat v_m(t)\,dt\\
+& \left(m-\frac 1 2\right)\log\frac{1+a}{1-a}-D_k +\int_{0}^\infty
\left[\frac{1}{t+D_k}-\frac{1}{t+2 A_m+\zeta}\right]\hat
v_m(t)\,dt\\
-&B_m
\end{split}
\end{equation}

Now,  we put
$$
C_m=\frac{2m-1}{1-a^2}.
$$
In this case the first line in \eqref{4.2} converges to $H_k(\zeta)$.
Since
\begin{equation}\label{}
\begin{split}
\lim_{m\to\infty}\int_{0}^\infty \left[\frac{1}{t+D_k}-\frac{1}{t+2
A_m+\zeta}\right] \left(\hat v_m(t)-\left(k+\frac 1 2\right)
\right)\,dt\\
=\int_{0}^\infty
\frac{\rho_k(t)-\left(k+\frac 1 2\right)}{t+D_k}dt
\end{split}
\end{equation}
and
\begin{equation}\label{}
\int_{0}^\infty \left[\frac{1}{t+D_k}-\frac{1}{t+2 A_m+\zeta}\right]
\,dt =\log\frac{2 A_m}{D_k}+\log\left(1+\frac{\zeta}{2 A_m}\right)
\end{equation}
we have from the second line in \eqref{4.2} that
\begin{equation}\label{}
\begin{split}
\lim_{m\to\infty}\left\{B_m-\left(m-\frac 1
2\right)\log\frac{1+a}{1-a}-\left(k+\frac 1 2\right)\log{2 A_m}\right\}\\
=-D_k-\left(k+\frac 1 2\right)\log D_k+\int_{0}^\infty
\frac{\rho_k(t)-\left(k+\frac 1
2\right)}{t+D_k}dt=-Y_k.
\end{split}
\end{equation}
By \eqref{f4} we get \eqref{4.1}.

Now, let us transform the convergence of conformal mappings in the
asymptotic for the extremal function. From \eqref{f.7} we have
\begin{equation}\label{f9}
    z=\sqrt{\frac{Z^2-A_m^2}{Z^2-C_m^2}}=\sqrt{\frac{2\zeta A_m+\zeta^2}{A_m^2-C_m^2+2\zeta A_m+\zeta^2}}
    \sim \sqrt{\frac{2a}{2m-1}}\sqrt{-\zeta}.
\end{equation}
From \eqref{f10} we get the following chain of equalities for the
depending variable
\begin{equation}\label{f11}
\begin{split}
    f(z;k,m;a)-1=&(-1)^{k+1}L_m^k\cos\phi_m\\
    =&(-1)^{k+1}L_m^k \cosh\Phi_m\\
    =&(-1)^{k+1}L_m^k \cosh(H^k_m+B_m).
    \end{split}
\end{equation}
Since $H^k_m\to H_k$, we have from \eqref{f3} and \eqref{3.2}
\begin{equation*}
    f(z;k,m;a)-1\sim (-1)^{k+1}\frac 1 \pi \int_0^\infty\frac
    {e^{\zeta-t} }{t-\zeta}\left(\frac{t}{-\zeta}\right)^{k-\frac 1
    2}dt.
\end{equation*}
Putting here $\zeta=-\lambda^2$ (see \eqref{f9}) and $t=\mu^2$, we
get \eqref{f12}.
\end{proof}

\section{Unweighted Extremal Polynomial via Conformal Mapping}
Let $P_m(z,p,a)$ be the best uniform (unweighted) approximation of
$|x|^p$ by polynomials of degree less or equal $2m$,  $2m>p$, on two
intervals $[-1,-a]\cup[a,1]$ and let $E=E_{2m}(p,a)$ be the
approximation error.

In this
section we prove
\begin{theorem}\label{t5.1} For a not even $p$
there is a curve $\gamma=\gamma_m(p,a)$ inside the half--strip
\begin{equation}\label{6.1.18}
\{w=u+iv: u\in (0,(m+1)\pi),\ v>0\}
\end{equation}
such that the extremal polynomial possesses the representation
\begin{equation}\label{fcm}
    P_m(z,p,a)=z^p+(-1)^{[p/2]}E\cos\phi_m(z,p,a)
\end{equation}
where $\phi_m(z,p,a)$ is the conformal map of the first quadrant
onto the region $\Omega_m(p,a)$ in the half strip \eqref{6.1.18}
bounded on the left by $\gamma_m(p,a)$. The conformal map is
normalized by $\phi_m(a,p,a)=0$, $\phi_m(1,p,a)=(m+1)\pi$ and
$\phi_m(\infty,p,a)=\infty$. Moreover, the curve $\gamma$ is the
image of the imaginary half--axis under this conformal map that
satisfies the following functional equation
\begin{equation}\label{6.2.18}
\gamma_m(p,a)=\{u+iv=\phi_m(iy,p,a):E\sin u(y)\sinh v(y)=\left|\sin
\frac {\pi p }2\right|y^p, y>0\}.
\end{equation}
\end{theorem}

\begin{proof}The proof contains two main ingredients: the Chebyshev
theorem and the argument principle. In addition to that we will show
some  particular fact related to  the shape of the extremal
polynomial. We prove that $P_m(0,p,a)>E$ for even $[p/2]$  and
$P_m(0,p,a)<-E$, when $[p/2]$ is odd.

 Due to the symmetry of  $P_m(x,p,a)$, we
can use the Chebyshev theorem with respect to the best approximation
of $(\sqrt{x})^p$ on $[a^2,1]$ by polynomials of degree $m$. It
gives us that $P_m(x,p,a)$ has $m+2$ points $\{x_j\}$  on the
interval $[a,1]$ where $P_m(x,p,a)-x^p$ alternates between $\pm E$
(the right half of the Chebyshev set in this case). Moreover,
$x_0=a$ and $x_{m+1}=1$. From this remark  we deduce that for
$|t|<1$ the equation
\begin{equation}\label{6.3.18}
P_m(x,p,a)-x^p=tE
\end{equation}
has  precisely $m+1$ zeros, say $\{x_j(t)\}$, on $(a,1)$. On the
other hand
$$
\{1,x^2,...,x^{2m}\}\cup\{x^p\}
$$
forms the so called Chebyshev system on $[0,\infty)$, see e.g.
\cite[Ch. II, Sect. 2]{KN}, and therefore \eqref{6.3.18} has  no
other solutions on the positive half
 axis ($m+1$ is the maximal possible number of roots for a
 generalized polynomial formed by a Chebyshev system
 of $m+2$ functions).

Using the argument principle we show that \eqref{6.3.18} has no
other solutions in the whole quarter--plane.

 Consider  the contour that runs on the positive real axis
till $x_j(t)-\epsilon$, then it goes around $x_j(t)$ on the
half--circle of the radius $\epsilon$ clockwise. After the last of
$x_j$'s we continue to go along the contour till the big positive
$R$. Next piece of the contour is a quarter--circle till imaginary
axis. Finally, from $iR$ we go back to the origin. On each
half--circle of the radius $\epsilon$ the argument of the function
changes by $-\pi$.  On the quarter circle it changes by about $\deg
P_m(z,p,a)\times\frac{\pi}2=m\pi$.

It remains to show that the change of the argument on the last piece
of the contour is about $\pi$. Then the whole change is
$-(m+1)\pi+m\pi+\pi=0$, and  since the function has no poles, it has
no zeros in the region.

Note that on the imaginary axis we have $\Re(P_m(iy,p,a)-(iy)^p)=
P_m(iy,a)-\cos{\frac {\pi p} 2}y^p$ and $\Im(P_m(iy,a)-(iy)^p)=
-\sin{\frac {\pi p } 2}y^p$. So the imaginary part increases with
$y$ for odd $\left[\frac p 2\right]$ and decreases when it is even.
Thus, it is enough  to show that the real part changes from a
certain negative value to $+\infty$ in the first case and, starting
from a positive value for $y=0$, it approaches to $-\infty$ as
$y\to\infty$ in the second case (recall that $2m>p$).

We give here a self--contained  proof of the above claim. For an
alternative proof see Remark   \ref{gk} below. Note that, if $a$ is
close to $1$, for $p=2k-1$ the shape of the extremal unweighted
polynomial is close to the shape of the extremal polynomial with the
weight $|x|^{2k-1}$, see Remark \ref{r1}. Consider, for example, the
first case, $k-1=\left[\frac p 2\right]$ is odd, then, due to
Theorem \ref{t3.1} and  the above remark,
\begin{equation}\label{f8}
    P_m(a,p,a)-a^p=-E\quad\text{and}\quad P_m(1,p,a)-1=(-1)^m E.
\end{equation}
Since $E_{2m}(p,a)\not=0$ for all $0<a<1$, $2k-2<p<2k$, no
bifurcation is possible and relations \eqref{f8} hold true for all
values of $(a,p)$ in the region. Since, moreover, \eqref{6.3.18} has
no solutions on $\bbR\setminus [a,1]$ we get the required behavior
of
 $P_m(z,p,a)-z^p$ as $z=iy$, $y\to 0$ and $y\to+\infty$, from its
 behavior on the real axis  $z=x$, as $x\to 0$ and $x\to+\infty$.

Thus $\arccos\frac{P_m(z,p,a)-z^p} E$ is well defined in the
quarter--plane. We finish the proof by  inspection of the boundary
correspondence.
\end{proof}

Note two facts: the curve \eqref{6.2.18} has the asymptote $u\to
\pi$, $v\to+\infty$ ($y\to+\infty$) and we have uniqueness of the
solution of the functional equation  \eqref{6.2.18} due to
uniqueness of the extremal polynomial.

\begin{remark}\label{gk} Recall Gantmacher--Krein's Theorem (see e.g.
\cite[Theorem 4.4, more specifically Corollary 4.4]{KSt}): the
number of distinct zeros on $(0,\infty)$ of any generalized
polynomial $\sum_{i=0}^n a_i x^{\alpha_i}$, where $\sum_{i=0}^n
a_i^2>0$ and $\alpha_0,\alpha_1...\alpha_n$ is an increasing
sequence of real numbers, is at most the number of sign changes in
the sequence $a_0,a_1,...,a_n$ after zero terms are discarded. Since
the "polynomial" $P_m(x,p,a)-x^p-tE$, $-1<t<1$, has the maximal
possible number of zeros in $(a,1)$ its coefficients sequence (in
the right order) has the maximal possible number of sign changes,
that is $m+1$. The coefficient before $x^p$ is negative, therefore
the signs of the zero coefficient, or $P_m(0,p,a)-tE$, and the last
one are $(-1)^{\left[\frac p 2\right]}$ and  $(-1)^{\left[\frac p
2\right]+m+1}$ respectively (compare \eqref{f8}).
\end{remark}

\section{And its Asymptotics}

\begin{theorem} For the approximation error $E_{2m}(p,a)$ the limit
\eqref{f1} exists. Moreover, uniformly on compact subsets of the
positive half--axis,
\begin{equation}\label{f13}
\begin{split}
\lim_{m\to\infty}&
    \left\{\left(\frac{m}{a}\right)^{\frac p 2}P_m\left(\sqrt{\frac{a}{m}}\lambda,p,a\right)\right\}\\
    =& \lambda^p+\frac {\sin\frac {\pi p} 2} \pi \int_0^\infty
    \mu^{p}e^{-(\lambda^2+\mu^2)}\frac
    {2\mu d\mu}
    {\lambda^2+\mu^2}.
    \end{split}
\end{equation}
\end{theorem}

\begin{proof}
First we present briefly the second  step of our method similar to
the proof of Theorem \ref{t4.1}. We use the representation
\eqref{fcm}. Then we use the symmetry principle and make convenient
changes of variable to have a conformal map $\Phi_m(Z)$ of the upper
plane in the region
$$
i(\Omega_m(p,a)\cup\overline{\Omega_m(p,a)})\cup(0,i\pi(m+1))
$$
with the boundary correspondence
$$
\Phi_m:(-C_m,-A_m,0,A_m,C_m)\to (-\infty,-B_m,0,B_m,\infty),
$$
here $A_m=aC_m$ and $C_m=\frac{2m}{1-a^2}$.

For $\Phi_m$ we have the integral representation
$$
\Phi_m(Z)=m\log\frac{1+\frac{Z}{C_m}}
{1-\frac{Z}{C_m}}+\int_{A_m}^\infty
\left[\frac{1}{X-Z}-\frac{1}{X+Z}\right]v_m(X)\,dX,
$$
where
\begin{equation}
v_m(X)=
\begin{cases}
\frac 1{\pi}\Im \Phi_m(X),& A_m\le X\le C_m \\
1,& X> C_m
\end{cases}
\end{equation}
and we put again
$$
H_m(\zeta)=\Phi_m(Z)-B_m,\quad Z=A_m+\zeta.
$$
Then
\begin{equation}\label{2f}
\begin{split}
H_m(\zeta)=&m\log\frac{1+a+\frac{\zeta}{C_m}}
{1-a-\frac{\zeta}{C_m}}+\int_{0}^\infty
\left[\frac{1}{t-\zeta}-\frac{1}{t+2A_m+\zeta}\right]\hat v_m(t)\,dt
-B_m,\\
\sim &m\log\frac{1+a} {1-a}+\zeta+\int_{0}^\infty
\left[\frac{1}{t-\zeta}-\frac{1}{t+2A_m+\zeta}\right]\hat v_m(t)\,dt
-B_m,
\end{split}
\end{equation}
In a usual way we write
\begin{equation}\label{3f}
\begin{split}
&\int_{0}^\infty
\left[\frac{1}{t-\zeta}-\frac{1}{t+2A_m+\zeta}\right]\hat
v_m(t)\,dt\\ & =\int_{0}^\infty
\left[\frac{1}{t-\zeta}-\frac{1}{t+2A_m+\zeta}\right](\hat
v_m(t)-\chi_{[1,\infty]}(t))\,dt\\& -\log(1-\zeta) +\log
(1+2A_m+\zeta)
\end{split}
\end{equation}

Since
$$
H_m(\zeta)\to w(\zeta)
$$
from \eqref{2f}, \eqref{3f} we have
\begin{equation}\label{fBm}
    B_m-m\log\frac{1+a} {1-a}-\log (2A_m)\to -c
\end{equation}
and
$$
w(\zeta)=\zeta-\log(-\zeta)+c+...
$$
as  $\zeta\to\infty$.

A bit new element:  rewrite the main equation \eqref{6.2.18} into
the form (the right hand side is not a constant any more)
\begin{equation}\label{1f}
E_{2m}\Im \cosh\Phi_m  =\left|\sin\frac{\pi p}2\right| y^p.
\end{equation}
 For $\zeta=\xi+i\eta$ we have
$$
y=\sqrt{\frac{(A_m+\xi)^2-A_m^2}{C_m^2-(A_m+\xi)^2}}\sim
\sqrt{\frac{2\xi a}{(1-a^2)C_m}} =\sqrt{\frac{\xi a}{m}}.
$$
Thus \eqref{1f} is of the form
\begin{equation}\label{main}
 \Lambda\Im e^w=\left|\sin\frac{\pi p}2\right|\xi^{\frac p 2}
\end{equation}
where
\begin{equation}\label{Ln}
    \Lambda=\lim_{m\to\infty}E_{2m}\left(\frac m a\right)^{\frac p 2}\frac 1 2e^{B_m}
\end{equation}

Finally, again, as the third step, we
  are looking for $w$ in the form
 \begin{equation}\label{form}
 w(\zeta)=\zeta+\log\left\{\frac
1{\pi}\int_0^\infty\frac{\tau(t)dt}{t-\zeta}\right\}.
 \end{equation}

Two small remarks on the normalization: due to $w(0)=0$ we have
\begin{equation}\label{norm}
    \frac 1{\pi}\int_0^\infty\tau(t)\frac{dt}{t}=1,
\end{equation}
and also
\begin{equation}\label{ec}
e^c=\frac 1{\pi}\int_0^\infty\tau(t){dt}.
\end{equation}

By  representation \eqref{form} the main equation \eqref{main} is
nothing but
\begin{equation}\label{lineq2}
    \left(e^\zeta\Im{\frac
1{\pi}\int_0^\infty\frac{\tau(t)dt}{t-\zeta}}\right)_{\zeta=\xi+i0}=
\frac{\left|\sin\frac{\pi p}2\right|}\Lambda \xi^{\frac p 2}.
\end{equation}
Thus
\begin{equation}\label{f5}
\tau(\xi)= \frac{\left|\sin\frac{\pi p}2\right|}\Lambda \xi^{\frac p
2}e^{-\xi}
\end{equation}
and basically we are done. By \eqref{ec}
\begin{equation*}
e^c=\frac{\left|\sin\frac{\pi p}2\right|}{\pi\Lambda}
\Gamma\left(\frac{p}{2} +1\right)=\frac
1{\Lambda\left|\Gamma\left(-\frac{p}{2}\right)\right|}.
\end{equation*}
The constant  $\Lambda$ is uniquely defined by \eqref{norm}. By
\eqref{Ln}, \eqref{fBm} we have
\begin{equation}\label{otvet}
    E_{2m}(p,a)\sim 2\Lambda\left(\frac{a}{m}\right)^{\frac p 2}\left(\frac{1-a}
    {1+a}\right)^m\frac {e^c}{2A_m}
    =
    \left(\frac{a}{m}\right)^{\frac p 2}\left(\frac{1-a}
    {1+a}\right)^m\frac {1-a^2}{2a m
\left|\Gamma\left(-\frac{p}{2}\right)\right|}
\end{equation}
and \eqref{f1} is proved.

The proof of \eqref{f13} is similar to the proof of \eqref{f12}.

\end{proof}

\section{Appendix 1}

Similar to \eqref{fcm} for extremal entire function $F$ and for
$a=0$ we write
\begin{equation}\label{al1}
    F(z)=z^p+(-1)^{[\frac p 2]}E\cos\phi(z),
\end{equation}
or for $z=-i\lambda$ and $\phi=-i\psi$
\begin{equation}\label{al2}
    F(i\lambda)=(-i\lambda)^p+(-1)^{[\frac p 2]}E\cosh\psi(\lambda).
\end{equation}
Now, $\psi(\lambda)$ is the conformal map of the upper half--plane
on the upper half-plane
\begin{equation}\label{al3}
\psi(\lambda)=\lambda+\log\left\{C+\frac 1
\pi\int\left(\frac{1}{\mu-\lambda}-\frac{\mu}{1+\mu^2}\right)\rho(\mu)d
\mu\right\}.
\end{equation}
Note that in this representation $\rho$ is {\it not} symmetric.

Put
$$
C+\frac 1
\pi\int\left(\frac{1}{\mu-\lambda}-\frac{\mu}{1+\mu^2}\right)\rho(\mu)d
\mu =-\tilde\rho+i\rho,
$$
where $\tilde\rho$ is "a kind of  Hilbert transform" of $\rho$. Then
we get from \eqref{al2}
\begin{equation}\label{al3}
    F(i\lambda)=e^{-i\frac{\pi p}2}\lambda^p+(-1)^{[\frac p 2]}\frac E 2
\left\{e^{\lambda}(-\tilde\rho+i\rho)+\frac{e^{-\lambda}}{-\tilde\rho+i\rho}\right\}.
\end{equation}
For  real $\lambda$'s the imaginary part of this expression gives us
\begin{equation}\label{al4}
    0=-\sin\frac{\pi p}2\lambda^p+(-1)^{[\frac p 2]}\frac E 2
\left\{e^{\lambda}\rho-\frac{e^{-\lambda}\rho}{\tilde\rho^2+\rho^2}\right\}.
\end{equation}

\begin{conjecture}\label{c}
For the Bernstein Problem, $a=0$, we conjecture that the extremal
entire function $F$ is of the form \eqref{al1}, where $\phi$ is the
conformal map  of the upper half plane onto the region in the upper
half plane above the curve
$$
\gamma=\{u+iv=\phi(x):\ x\in\bbR\}
$$
such that
\begin{equation}\label{trigeq0}
    L\sin v(x)\sinh u(x)=x,\quad x\in \bbR
\end{equation}
and normalized by $\phi(0)=0$, $\phi(z)\sim z$, $z\to\infty$.
\end{conjecture}

Let us rewrite the above equation in terms of the unknown function,
say $\rho$ and its Hilbert transform $\tilde\rho$. We use the
integral representation
$$
\phi(z)=z+\frac 1 \pi \int_0^\infty
\left[\frac{1}{x-z}-\frac{1}{x+z}\right]v(x)\,dx.
$$
The curve  has the asymptote $v\to \pi$, $x\to\infty$. We define
$\rho:=\pi-v$ to write
$$
\phi(z)=z+i\pi-\frac 1 \pi \int_0^\infty \frac{\rho\,dx}{x-z}.
$$
Finally since
$$
\frac 1 \pi \int_0^\infty \frac{\rho\,dx}{x-z}=-\tilde\rho+i\rho,
$$
we get $u(x)=\tilde\rho(x)+x$ and $v(x)=\pi-\rho(x)$. Thus equation
\eqref{trigeq0} leads to
\begin{equation}\label{trigeq}
    L\sin\rho(x)\sinh(\tilde\rho(x)+x)=x.
\end{equation}

\section{Appendix 2}
From \cite{Akha}, problem 42:
\begin{equation}\label{6.1}
E_l\left[\frac{1}{(b+x)^s}\right] \sim \frac{l^{s-1}}{|\Gamma(s)|}
\frac{(b-\sqrt{b^2-1})^l} {(b^2-1)^{\frac{s+1}{2}}}\quad (b>1, \
s\not= 0),
\end{equation}
where $E_l[f(x)]$ is the error of the approximation of the function
$f(x)$ on the interval $[-1,1]$ by polynomials of degree not more
than $l$.

We change the variable
$$
y=\frac{b+x}{b+1}
$$
and put $a^2=\frac{b-1}{b+1}$. Then we have
$$
\inf_{P:\deg P\le l}\max_{y\in [a^2,1]}|y^{-s}-P(y)| =(1+b)^s
E_l\left[\frac{1}{(b+x)^s}\right].
$$
That is
\begin{equation}\label{6.2}
E_{2l}(-2s,a)=(1+b)^s E_l\left[\frac{1}{(b+x)^s}\right].
\end{equation}

Note that
$$
b=\frac{1+a^2}{1-a^2},\quad b^2-1=\frac{4a^2}{(1-a^2)^2},
$$
and therefore
$$
\sqrt{b^2-1}=\frac{2a}{1-a^2}, \quad b-\sqrt{b^2-1}=\frac{1-a}{1+a}.
$$
Thus from \eqref{6.1} and \eqref{6.2} we get
\begin{equation*}
\begin{split}
E_{2l}(-2s,a)\sim & \left(\frac{2} {1-a^2}\right)^{s}
\frac{l^{s-1}}{|\Gamma(s)|} \left(\frac{1-a} {1+a}\right)^l
\left(\frac{1-a^2}
{2a}\right)^{s+1}\\
= & a^{-s} \frac{l^{s-1}}{|\Gamma(s)|} \left(\frac{1-a}
{1+a}\right)^l \left( \frac{1-a^2} {2a}
\right)\\
=& a^{-s-1} \frac{l^{s-1}}{|\Gamma(s)|} \left(\frac{1-a}
{1+a}\right)^{l+1} \frac{(1+a)^2} 2.
\end{split}
\end{equation*}

\bibliographystyle{amsplain}

\noindent
Addresses: \\[.1cm]
Fedor Nazarov \\
Department of Mathematics \\
Michigan State University \\
East Lansing, Michigan 48824, USA \\
fedja@math.msu.edu \\

\noindent
Franz Peherstorfer \\
Abteilung f\"ur Dynamische Systeme \\
und Approximationstheorie\\
Universit\"at Linz\\
4040 Linz, Austria\\
{Franz.Peherstorfer@jku.at}\\

\noindent
Alexander Volberg\\
Department of Mathematics \\
Michigan State University \\
East Lansing, Michigan 48824, USA \\
volberg@math.msu.edu\\

\noindent
Peter Yuditskii\\
Abteilung f\"ur Dynamische Systeme \\
und Approximationstheorie\\
Universit\"at Linz\\
4040 Linz, Austria\\
 Petro.Yudytskiy@jku.at

\end{document}